\documentclass[10pt]{article}
\usepackage[english,francais]{babel}
\usepackage{fullpage}
\advance\hoffset by -3mm  
\advance\voffset by  8mm  

\usepackage{graphicx}

\usepackage{amssymb}

\def\ha{{\widehat{a}}}
\def\hA{{\widehat{A}}}
\def\hb{{\widehat{b}}}
\def\hh{{\widehat h}}

\def\cO{{\mathcal O}}

\usepackage[english,francais]{babel}

\newtheorem{e-proposition}[theorem]{Proposition}

\newtheorem{e-definition}[theorem]{Definition\rm}

\newtheorem{exemple}{\it Exemple\/}
\renewcommand{\theequation}{\arabic{equation}}
\def\og{\leavevmode\raise.3ex\hbox{$\scriptscriptstyle\langle\!\langle$~}}
\def\fg{\leavevmode\raise.3ex\hbox{~$\!\scriptscriptstyle\,\rangle\!\rangle$}}

\begin{document}

\selectlanguage{francais} \title{Une m{\'e}thode d'identification
  pour un syst{\`e}me lin{\'e}aire {\`a} retards\\
\ \\
\selectlanguage{english}
A method of identification for a linear delays system.}



\author{
 Fran{\c{c}}ois Ollivier\\
 ALIEN, INRIA Futurs \&\ LIX, UMR CNRS 7161\\
 {\'E}cole polytechnique\\ 
 91128 Palaiseau CEDEX, France\\
{\tt francois.ollivier@lix.polytechnique.fr}\\
\ \\
Sa{\"\i}d Moutaouakil et Brahim Sadik\\
D{\'e}partement de Math{\'e}matiques\\
 Facult{\'e} des Sciences Semlalia\\
 B.P. 2390, Avenue Safi, Marrakech, Maroc\\
{\tt $\{$sadik,s.moutaouakil$\}$@ucam.ac.ma}
 }

\date{\today}

\maketitle


\selectlanguage{english}
\begin{abstract} 

We provide a class of methods for the identification of a
  linear system with delay of the shape $x^{(n)}(t) =
\sum_{i=0}^{n-1}a_{i}x^{(i)}(t) + bu(t-h)$. They allow the
simultaneous identification of the parameters and delay, the
observation of $x$ and its derivatives, knowing only generic input $u$
and output $x$. They are robust to the noise. Used in continuous time,
they may follow the evolution of slowly varying parameters and
noise. The method may be generalized to systems with two or more
delays, e.g. $x''(t)+ax(t-h_{1})=bu(t-h_{2})$.

\vskip 0.5\baselineskip

\selectlanguage{francais}
\noindent{\bf R\'esum\'e} On introduit une classe de m{\'e}thodes
d'identification et d'observation pour un syst{\`e}me lin{\'e}aire
{\`a} retard de la forme $x^{(n)}(t) = \sum_{i=0}^{n-1}a_{i}x^{(i)}(t)
+ bu(t-h)$. Celles-ci permettent l'identification simultan{\'e}e du
retard et des param{\`e}tres, l'observation de $x$ et de ses
d{\'e}riv{\'e}es et ne supposent que la connaissance de la sortie $x$
et de l'entr{\'e}e $u$, sans autre hypoth{\`e}se que leur
g{\'e}n{\'e}ricit{\'e}. Elles se r{\'e}v{\`e}lent assez robustes au
bruit. Utilis{\'e}es en temps continu, elles permettent de suivre
l'{\'e}volution de param{\`e}tres ou de retards lentement variables.
La m{\'e}thode peut {\^e}tre g{\'e}n{\'e}ralis{\'e}e {\`a} des
syst{\`e}mes avec plusieurs retards,
e.g. $x''(t)+ax(t-h_{1})=bu(t-h_{2})$.

\noindent

\end{abstract}

\selectlanguage{english}
\section*{Abridged English version}

Time delay systems are known to take an important place in many fields
of application. We refer to \cite{Richard2003} and the references
therein  for more details on the subject. We propose a method for
parameters and delay identification, related to the method introduced by Fliess
and Sira-Ram{\`\i}rez \cite{Fliess-Sira2003a,Fliess-Sira2003b} (see also
\cite{Belkoura-Richard-Fliess2006}).
\medskip

We write the system $\sum_{i=0}^{n} a_{i}x^{(i)}(t) + b u(t)=0$, with
$a_{n}=-1$.  The output is $y=x+w$ where $w$ stands for the noise
which is such that $\left(\int_{T_{1}}^{T_{2}}
w(\tau)d\tau\right)/(T_{1}-T_{2})$ goes to $0$ when $T_{2}-T_{1}$ goes
to infinity. The main idea is to use a family of functions $f_{i}$
such that $f_{j}^{(k)}(T_{1})=f_{j}^{(k)}(T_{2})=0$ for $k<n$.  Let
$I_{x,f}:=\int_{T_{1}}^{T_{2}} f(\tau)x(\tau)d\tau$. Integrating by
parts, we get $I_{x^{(i)},f_{j}}=(-1)^{i}I_{x,f_{j}^{(i)}}$.  So we
may estimate the values of the coefficients $a_{i}$ and $b$ by solving
the system $ \sum_{i=0}^{n} \left((-1)^{i}a_{i}
I_{y,f_{j}^{(i)}}\right) + b I_{u,f_{j}}= 0 $ for $j=1,\ldots,m$, by
the mean squares method. We can estimate $x$ and its derivatives in
the same way, using functions $g_{j}$ $0\le j< n$ such that
$g_{j}^{(k)}(T_{1})=0$ for $0\le k< n$ and $g_{j}^{(k)}(T_{2})=0$ for
$0\le k< j$, with $g_{j}^{(j)}(T_{2})\neq0$.

In practice we have used $f_{j}(t)=
(T_{2}-t)^{n+j}e^{-\lambda(T_{2}-t)}$, the integration being done
between $-\infty$ and the current time. A good approximation of the
integrals is obtained by integrating the system $J_{x,0}'=x - \lambda
J_{x,0}$ and $J_{x,j}'=jJ_{x,j-1} - \lambda J_{x,j}$ if $j>0$, with
initial conditions $J_{j}(0)=0$: $J_{j}(t)$ tends quickly to
$I_{x,f_{j-n+1}}$ for $\lambda$ great enough.  Numerical simulations
are given with example 1.
\medskip

In section~3, we consider a delay system $\sum_{i=0}^{n}
a_{i}x^{(i)}(t) + b u(t-h)$, with $a_{n}=-1$. We use the notation
$I_{x,f,T_{1},T_{2}}=\int_{T_{1}}^{T_{2}}f((\tau-T_{1})/(T_{2}-T_{1}))
x(\tau)d\tau$, where  $f$ is such that $f^{(j)}(0)=f^{(j)}(1)=0$, for
$j<m$ with $m\ge n$. Let $E_{(T_{1},T_{2}),x,u}$ denote the equation
$$
\sum_{i=0}^{n} (T_{1}-T_{2})^{-i}a_{i}I_{x,f^{(i)},T_{1},T_{2}} 
+ \sum_{\ell=0}^{k}b_\ell (T_{2}-T_{1})^{-\ell}
I_{u,f^{(\ell)},T_{1},T_{2}} + O(h^{\min(k+1,m)})
=0,
$$ where the term $O(h^{\min(k+1,m)})$ is neglected.  Solving it by the mean
squares method for generic $x$ and $u$ and a generic set $S$ of couples
$(T_{1},T_{2})$, we get approximations $\ha_{i}$ of the coefficients
and an approximation of delay equal to $\hh=\hb_{1}/\hb_{0}$. We can
then replace $u$ by $u_{\hh}(t)=u(t-\hh)$ in $E_{(T_{1},T_{2}),x,u}$
in order to get the improved approximation and iterate the process.
Better precision could be achieved if $f$ is such that
$f(t-h)=\sum_{k=0}^{p} c_{k}(h)f^{(k)}(t)$, for example with
$f(t)=sin^{m}(\pi t)$.

In example~2, we use $f=\sin^{2}$ on $[0, \pi]$ and investigate the precision of
the evaluation depending on the size of the noise. See {\it Table~0}
and {\it Table~1}, which shows how our method may be improved by using
a nonlinear least quare method, initialized with our data.

In example~3 the method is generalized to the two delays system
$x''(t)+ax(t-h_{1})=bu(t-h_{2})$. The standard deviation of our
computed values are given in {\it table~2} for different noise level,
stared values being improved by using a nonlinear least square method
initialized with our results.
\medskip

In section~4, we adapt the method of section~1 to the delay situation,
solving the system $$ \sum_{i=0}^{n} \left((-1)^{i}a_{i}
I_{y,f_{j}^{(i)}}\right) + b_{0}I_{u_{\hh},f_{j}}
+b_{1}I_{u_{hh},f_{j}} +O(h^{2})= 0,
$$ where $u_{\hh}(t)=u(t-\hh)$ for $|h-\hh|\ll 1$. The delay
evaluation is then $\hh+\hb_{1}/\hb_{0}$. For $t>T_{0}$, when
$\hb_{0}$ does not vanish any more and the evaluation
$\hh(0)+\hb_{1}/\hb_{0}$ is assumed to be close enough, we take
$\hh'=\lambda_{h}\hb_{1}/\hb_{0}$, so that $\hh$ will converge to $h$.

Example 3 shows a simulation with $h=0.5$. Greater delays could be
considered by changing the time scale. Example 4 considers the case of
a slowly varying delay. Scilab simulation files are available at url
\cite{Simus}.

\selectlanguage{francais}
\section{Introduction}
\label{intro}

Les syst{\`e}mes {\`a} retard sont l'objet de recherches actives, en
raison de leur importance dans de nombreux champs d'applications. Nous
renvoyons {\`a} l'article de synth{\`e}se \cite{Richard2003} et aux
r{\'e}f{\'e}rences incluses pour plus de d{\'e}tails. L'identification
simultan{\'e}e du retard et des param{\`e}tres demeure un probl{\`e}me
difficile, mais crucial pour la mise en {\oe}uvre de m{\'e}thodes de
contr{\^o}le appropri{\'e}es.

Nous proposons une m{\'e}thode qui peut s'interpr{\'e}ter comme une
variante de la m{\'e}thode d'identification et d'observation due {\`a} Michel
Fliess et Hebertt Sira-Ram{\`\i}rez
\cite{Fliess-Sira2003a,Fliess-Sira2003b}, et qui est susceptible de
s'adapter {\`a} des syst{\`e}mes {\`a} retard.  L'article
\cite{Belkoura-Richard-Fliess2006} fournit une application plus
directe de cette approche. On pourra consulter, par exemple,
\cite{Tuch-Feuer-Palmor1994,Verduyn-Lunel1997,Tan-Wang-Lee,Mahji-Atherton1999,Kolmanovskii-Myshkis1999,Diop-Kolmanovsky-Moraal-vanNieuwstadt2001,Orlov-Belkoura-Richard-Dambrine2002}
pour un aper{\c{c}}u d'autres approches.

L'introduction d'un facteur d'oubli permet l'identification et
l'observation en continu, ce qui est d{\'e}crit dans la section 2 pour
des syst{\`e}mes sans retard et dans la section 4 pour des
syst{\`e}mes {\`a} retard.

La section 3 d{\'e}crit une m{\'e}thode d'identification {\it a
posteriori}. Une simulation permet d'{\'e}valuer sa pr{\'e}cision en
fonction du niveau de bruit. Nous avons utilis{\'e} {\`a} titre
d'exemple des syst{\`e}mes d'ordre $2$, qui correspondent {\`a} une
classe de mod{\`e}les tr{\`e}s r{\'e}pandus (voir
\cite{Fliess-Mounier2000}). Les calculs num{\'e}riques ont {\'e}t{\'e}
effectu{\'e}s sur Dell Latitute D400 avec Scilab~3.0, {\`a} partir de
formules litt{\'e}rales calcul{\'e}es en Maple~9.5.

\section{Identification d'un syst{\`e}me sans retard}

Pour all{\'e}ger les formules, nous {\'e}crirons le syst{\`e}me sous la forme:
$\sum_{i=0}^{n} a_{i}x^{(i)}(t) + b u(t)=0$, avec $a_{n}=-1$.  La
sortie est $y=x+w$, o{\`u} $w$ d{\'e}signe le bruit, que nous
supposons tel que $\left(\int_{T_{1}}^{T_{2}}
w(\tau)d\tau\right)/(T_{1}-T_{2})$ tend vers $0$ quand $T_{2}-T_{1}$
tend vers l'infini. L'id{\'e}e de base est
d'utiliser une famille de fonctions $f_{j}$, $j=1,\ldots,m$ telles que
$f_{j}^{(k)}(T_{1})=f_{j}^{(k)}(T_{2})=0$ pour $k<n$.

On pose $I_{x,f}:=\int_{T_{1}}^{T_{2}}
f(\tau)x(\tau)d\tau$. En int{\'e}grant par parties, on aura
$I_{x^{(i)},f_{j}^{(k)}}= -I_{x^{(i-1)},f_{j}^{(k+1)}}$, si $k<n$, et donc
$I_{x^{(i)},f_{j}}=(-1)^{i}I_{x,f_{j}^{(i)}}$. On peut donc calculer la valeur
des coefficients $a_{i}$ et $b$ en r{\'e}solvant le syst{\`e}me
$$ \sum_{i=0}^{n} \left((-1)^{i}a_{i} I_{y,f_{j}^{(i)}}\right) + b
I_{u,f_{j}}= 0
$$ pour $j=1,\ldots,m$, si $m\ge n$. Pour $m>n$, ce syst{\`e}me sera
r{\'e}solu par la m{\'e}thode des moindres carr{\'e}s.

Si l'on se donne  une famille $g_{j}$, $0\le j< n$ telle que
$g_{j}^{(k)}(T_{1})=0$ pour  $0\le k< n$ et $g_{j}^{(k)}(T_{2})=0$
pour  $0\le k< j$, avec $g_{j}^{(j)}(T_{2})\neq0$, on obtient alors
$I_{x^{(i)},g_{j}}=(-1)^{i}I_{x,g_{j}^{(i)}}
+\sum_{\ell=0}^{i-j}(-1)^{i-\ell-1}g_{j}^{(i-\ell-1)}(T_{2})x^{(\ell)}(T_{2})$.
On estime $x(T_{2}), \ldots, x^{(n-1)}(T_{2})$ en
r{\'e}solvant le syst{\`e}me
$$ \sum_{i=0}^{n}
\left(\sum_{\ell=0}^{i-j}
\left((-1)^{i-\ell-1}g_{j}^{(i-\ell-1)}(T_{2})x(T_{2})^{(\ell)}\right)
+(-1)^{i}a_{i}I_{y,g_{j}^{(i)}}\right) + b I_{u,g_{j}}=0,
$$ pour $j=0, \ldots,n-1$. Ces {\'e}quations se d{\'e}duisent
imm{\'e}diatement de
$\sum_{i=1}^{n}a_{i}I_{x^{(i)},g_{j}}+bI_{u,g_{j}}=0$.
\medskip

Par exemple, on peut prendre $f_{j}(t)=
(T_{2}-t)^{n+j}e^{-\lambda(T_{2}-t)}$, et
$g_{j}=(T_{2}-t)^{j}e^{-\lambda(T_{2}-t)}$. Ces fonctions satisfont
les hypoth{\`e}ses pour $T_{1}=-\infty$. Si l'on pose $J_{x,0}'=x -
\lambda J_{x,0}$ et $J_{x,j}'=jJ_{x,j-1} - \lambda J_{x,j}$ si $j>0$, avec
les conditions initiales $J_{j}(0)=0$, $J_{j}(T_{2})$ tend vers
$I_{x,g_{j}}$ pour $j=1,\ldots,n-1$ et vers $I_{x,f_{j-n+1}}$ pour
$j\ge n$, quand $T_{2}$ tend vers l'infini. La convergence est rapide
pour $\lambda$ suffisament grand. On a $g_{j}'=\lambda g_{j}-jg_{j-1}$
et une relation similaire pour les $f_{j}$. On a donc
$I_{x^{(i)},f_{j}}=\sum_{k=j-1}^{n-1}
M_{i,j,k}I_{x,g_{k}}+\sum_{k=1}^{j} M_{i,n+j-1,k}I_{x,f_{k}} $, o{\`u}
les matrices $M_{i,j,k}$ satisfont des relations de r{\'e}currence
simples, ce qui facilite les calculs.
On obtient en un temps $t$ g{\'e}n{\'e}rique et suffisament grand une
bonne approximation des param{\`e}tres $a_{i}$ et $b$ en r{\'e}solvant
le syst{\`e}me
$$ \sum_{i=0}^{n} a_{i} \left(\sum_{k=j-1}^{n-1}
M_{i,j,k}I_{y,g_{k}}+\sum_{k=1}^{j} M_{i,n+j-1,k}I_{y,f_{k}}\right) +
b I_{u,f_{j}}= 0,
$$ que nous noterons $CA=D$, o{\`u} $A=(a_{1}, \ldots, a_{n-1},
b)^{t}$.  Comme il se peut que pour certaines valeurs le rang du
syst{\`e}me ne s'annule, il est pr{\'e}f{\'e}rable de ne pas le
r{\'e}soudre directement mais d'en chercher la solution par une
m{\'e}thode de gradient, en se donnant de nouvelles variables $\ha_{i}$ et
$\hb$, et en int{\'e}grant le syst{\`e}me $\hA' = -\Lambda C^{\rm t}(CA -
D)$.  Ce choix permet de conserver des valeurs pr{\'e}cises lorsque
$C$ devient mal conditionn{\'e}e et de r{\'e}duire encore l'influence
du bruit.  Les valeurs de $\lambda$ et de $\Lambda$ r{\'e}sultent d'un
compromis entre la vitesse de convergence souhait{\'e}e et la
pr{\'e}cision recherch{\'e}e.
\smallskip
\begin{exemple}
On a pris $x$ solution de l'{\'e}quation $x''=a_{1}x'+a_{0}x+bu$,
o{\`u} $u=60\cos(1,23t+1,3\sin(t)-0,7\cos(0,5t))$, avec les conditions
initiales $x(0)=20$ et $x'(0)=0,3$. On a choisi un bruit gaussien avec
un {\'e}cart type de $5$. La sortie est {\'e}chantillonn{\'e}e {\`a}
$100$~Hz. On a pris $b=2$, $a_{0}=-0,35$, $a_{1}(t)=-1,2$ si $t<30$ et
$a_{1}=-1,2-0,02(t-30)$ sinon. Les courbes de la
figure~\ref{courbes-exemple1} r{\'e}sument les r{\'e}sultats obtenus
pour $\lambda=1$ et $\Lambda=10^{-3}$.
\begin{figure}[!ht]\label{courbes-exemple1}
\begin{center}\vbox to 0.4truecm{\vfill}
\hbox to \hsize{\includegraphics[scale=0.35]{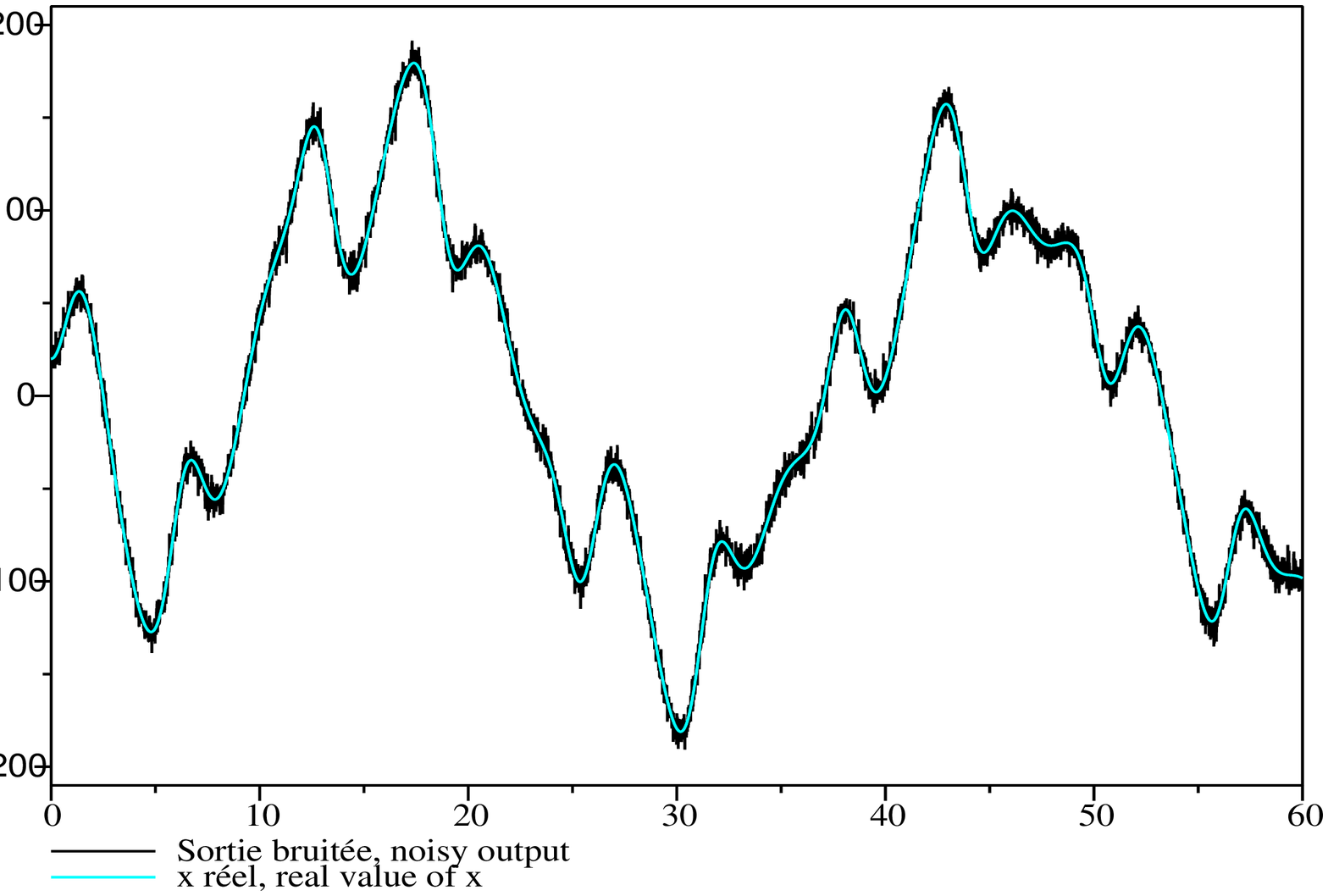}\hss
\includegraphics[scale=0.35]{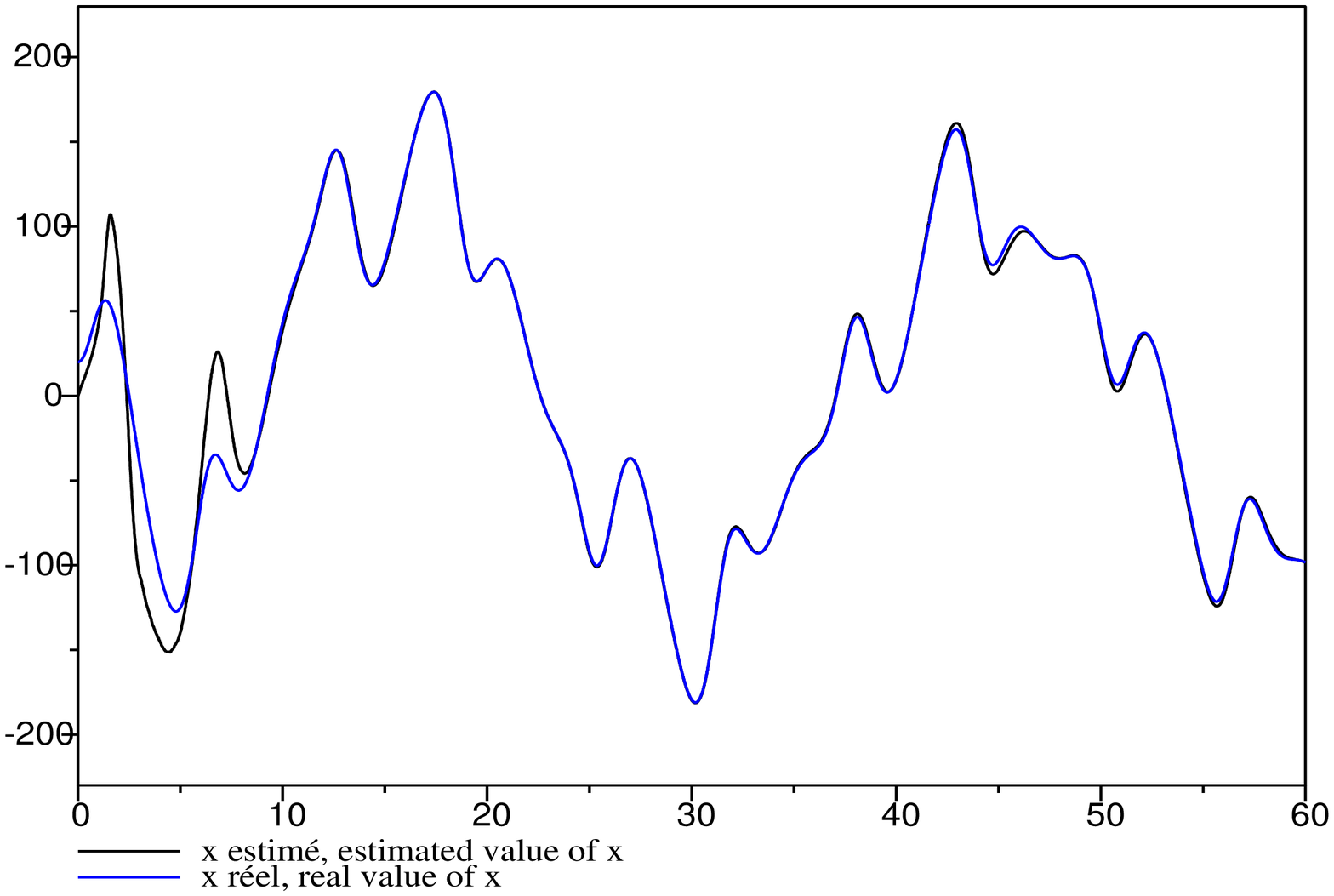}}
\medskip
\hbox to \hsize{\includegraphics[scale=0.35]{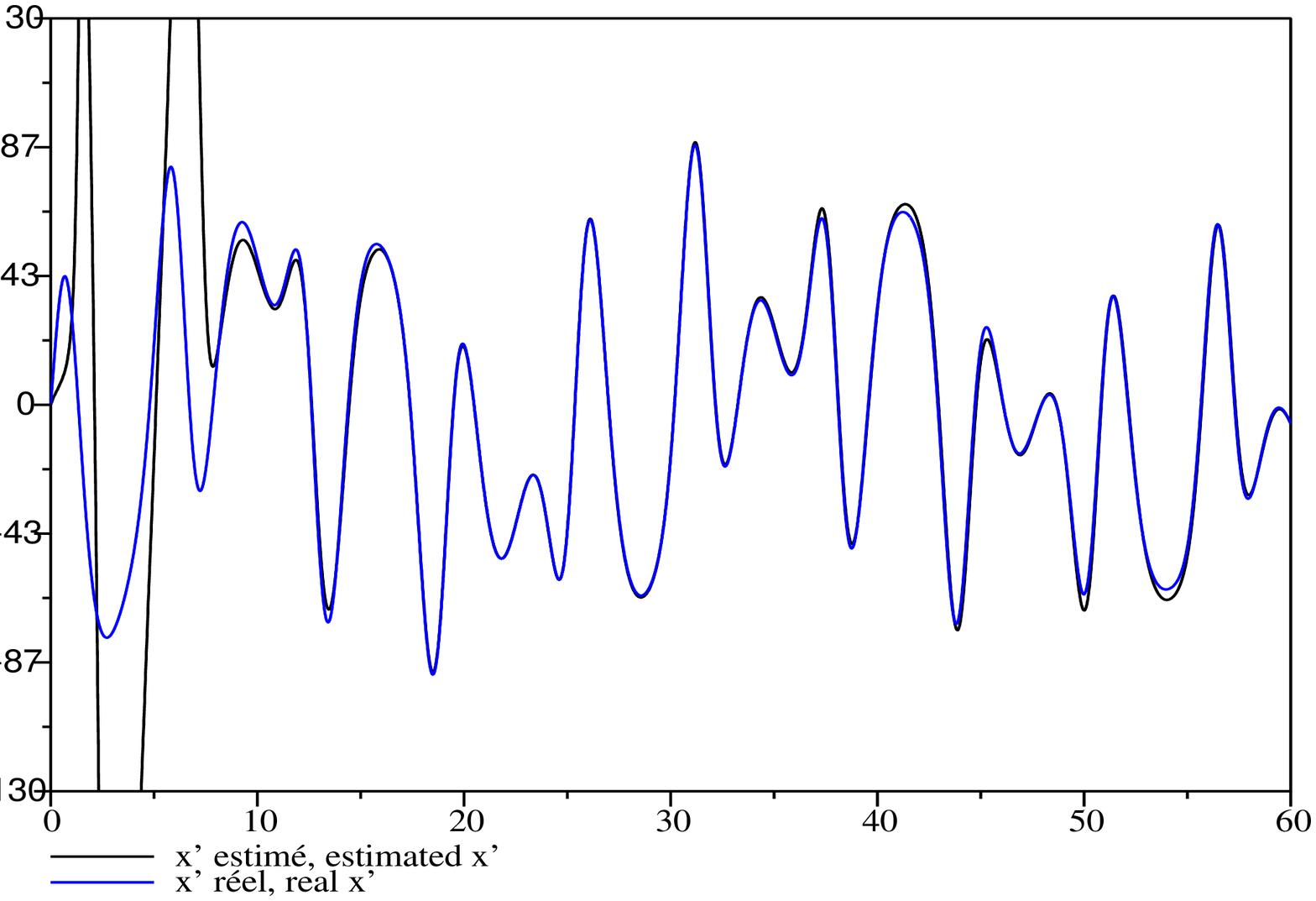}\hss
\includegraphics[scale=0.35]{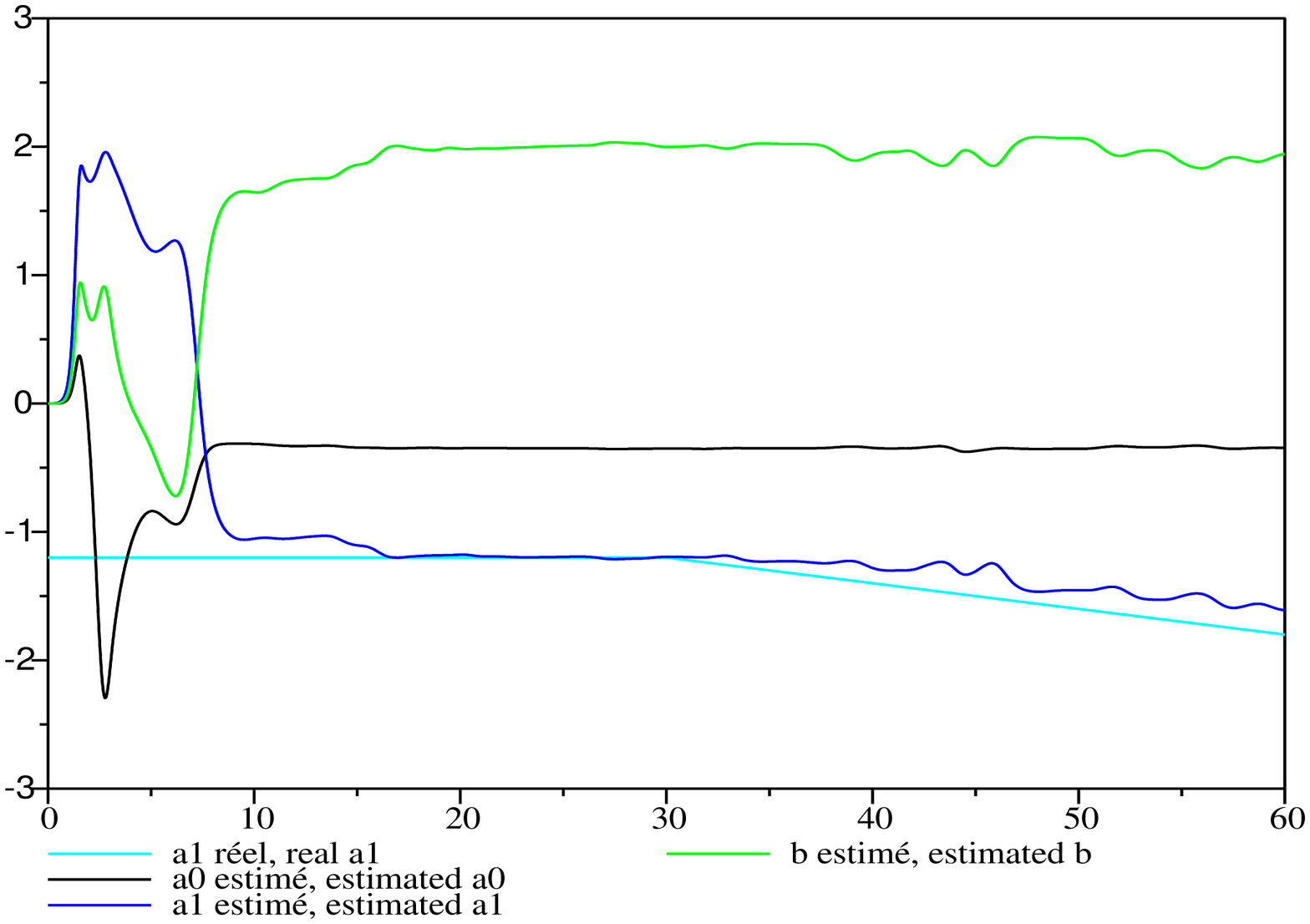}}
\caption{Estimations des coefficients de $x$ et de $x'$, {\it
Estimations of the coefficients, $x$ and $x'$.}}
\end{center}
\end{figure}
On voit qu'apr{\`e}s une phase de convergence,
l'approximation obtenue est excellente. La pr{\'e}cision de
l'{\'e}valuation des coefficients est naturellement moins bonne quand
$a_{1}$ commence {\`a} varier. Toutefois, les {\'e}valuations de $x$
et $x'$ demeurent assez pr{\'e}cises.
\end{exemple}

\vskip-0.3truecm
\section{Identification avec retards}

On consid{\`e}re maintenant un syst{\`e}me {\`a} retard {\'e}crit sous la forme
$\sum_{i=0}^{n} a_{i}x^{(i)}(t) + b u(t-h)$, avec $a_{n}=-1$.
Soit $f$ une fonction telle que $f^{(j)}(0)=f^{(j)}(1)=0$, pour
$j<m$ avec $m\ge n$. On pose
$I_{x,f,T_{1},T_{2}}=\int_{T_{1}}^{T_{2}}f((\tau-T_{1})/(T_{2}-T_{1}))
x(\tau)d\tau$.  
En int{\'e}grant par parties, on d{\'e}duit de l'{\'e}quation du syst{\`e}me
$$
\sum_{i=0}^{n} (T_{1}-T_{2})^{-i}a_{i}I_{x,f^{(i)},T_{1},T_{2}} 
+ b \left(\sum_{\ell=0}^{k}(T_{2}-T_{1})^{-\ell}
I_{u,f^{(\ell)},T_{1},T_{2}}h^{\ell}\right) + O(h^{\min(k+1,m)})
=0.
$$ On d{\'e}signera par $E_{(T_{1},T_{2}),x,u}$ cette {\'e}quation,
o{\`u} le terme $O(h^{\min(k+1,m)})$ a {\'e}t{\'e} n{\'e}glig{\'e} et
chaque puissance $h^{\ell}$ remplac{\'e}e par $b_{\ell}/b$ afin de se
ramener {\`a} une {\'e}quation lin{\'e}aire.  Pour $x$ et $u$
g{\'e}n{\'e}riques, on r{\'e}soud par moindres carr{\'e}s le
syst{\`e}me $E_{s,y,u}$, $s\in S$ en les inconnues $a_{i}$ et
$b_{\ell}$, o{\`u} $S$ est un ensemble g{\'e}n{\'e}rique
d'intervalles d'un cardinal suffisant. On obtient des approximations
$\ha_{i}$ et $\hb_{0}$ des coefficients $a_{i}$ et $b$ et une
approximation du retard {\'e}gale {\`a} $\hb_{1}/\hb_{0}$. La
pr{\'e}cision du r{\'e}sultat sera d'autant meilleure que le rapport
$h/(T_{2}-T_{1})$ est petit pour chaque couple de $S$.

Une fois obtenue une approximation $h_{0}$ du retard, on am{\'e}liore
la pr{\'e}cision du r{\'e}sultat en r{\'e}solvant le syst{\`e}me
$E_{s,y,u_{h_{0}}}$, $s\in S$, o{\`u} $u_{h_{0}}(t)=u(t-h_{0})$, ce
qui donnera la nouvelle approximation $h_{0}+\hb_{1}/\hb_{0}$.
On it{\`e}re le processus.

Pour de meilleurs r{\'e}sultats, on peut {\'e}galement choisir une fonction
$f$ telle que $f(t+h)=\sum_{k=0}^{p} a_{k}(h) f^{(k)}$, par exemple
$f(t)=t^{m}(1-t)^{m}$ ou $\sin^{m}(\pi t)$.
\medskip

\begin{exemple}
Nous avons fait des simulations num{\'e}riques avec $f(t)=\sin^{2}(t)$
sur l'intervalle $[0,\pi]$ pour un syst{\`e}me d'ordre
$2$. L'{\'e}quation $E_{(T_{1},T_{2}),y,u}$ s'{\'e}crit alors $
(T_{1}-T_{2})^{-1}a_{1}I_{x,f',T_{1},T_{2}}$
$+a_{0}I_{x,f,T_{1},T_{2}}+b_{0}I_{u,f,T_{1},T_{2}}$
$-b_{1}I_{u,cos(2t),T_{1},T_{2}}$
$+b_{2}I_{u,sin(2t),T_{1},T_{2}}+0(h^{2})$
$=(T_{1}-T_{2})^{-2}I_{x,f'',T_{1},T_{2}}$ et l'on obtient une
approximation de $h$ {\'e}gale {\`a} $h_{0}+{\rm
acos}(\hb_{2}/\hb_{0})$ {\`a} partir d'un syst{\`e}me
$E_{s,u_{h_{0}}}$, $s\in S$, o{\`u} $h_{0}$ est une estimation
pr{\'e}alable du retard.

Les valeurs des coefficients sont les m{\^e}mes que celles de l'exemple
$1$. On a choisi la commande
$u=60\cos(1,23t+0,33*\sin(t)-0,47\cos(0,5t))$ et un retard $h=4$. Le
bruit est un bruit gaussien et la fr{\'e}quence d'{\'e}chantillonnage est de
$500$~Hz. On a pris pour $S$ les couples $(10k,10k+15)$, $1\le k\le
9$. On part de $h_{0}=0$ et l'on it{\`e}re les calculs une cinquantaine de
fois (en un temps de l'ordre de quelques secondes sur un Dell Latitude
D400), ce qui est suffisant pour atteindre un point fixe {\`a} la
pr{\'e}cision num{\'e}rique pr{\`e}s.

Le tableau ci-dessous, qui donne la moyenne (M.) et l'{\'e}cart type
(\'E.) de s{\'e}ries de cent estimations des coefficients et du
retard r{\'e}alis{\'e}es pour des bruits gaussiens d'{\'e}carts types
croissant, montre que l'on obtient des estimations cr{\'e}dibles,
m{\^e}me pour un niveau de bruit notable.
\smallskip

{\em \scriptsize
\hbox to \hsize
{\hss\begin{tabular}{|c|c|c|c|c|c|c|c|c|}
\hline
\'E.~bruit & M.~$a_{1}$ & \'E.~$a_{1}$ & M.~$a_{0}$ & \'E.~$a_{0}$ & M.~$b$ & \'E.~$b$ & M.~$h$ & \'E.~$h$\\
\hline
1 &  -1,20 & 0,020  & -0,35 & 0,007  & 2,00  & 0,04 &  4,00 & 0,015\\
\hline
2 & -1,20  & 0,028  & -0,35 & 0,010  & 2,00  & 0,05 &  4,00 & 0,023\\
\hline
5 & -1,17  & 0,084  & -0.34 & 0,029  & 1,94  & 0,17 &  3,97  & 0,065\\
\hline
10 & -1,12 & 0,15   & -0,32 & 0,055  & 1,84  & 0,33 &  3,91  & 0,132\\
\hline
\end{tabular}\hss}
\hbox to \hsize{\hss Tableau 0. Moyenne et {\'e}cart type des valeurs
calcul{\'e}es\hss} 
\hbox to \hsize{\hss \em Table 0. Mean and standard
deviations {\em (\'E.)} of computed values\hss}}

\end{exemple}

Nous avons compar{\'e} nos r{\'e}sultats avec une m{\'e}thode souvent
utilis{\'e}e en pratique qui consiste {\`a} {\'e}valuer les
coefficients par une m{\'e}thode de moindres carr{\'e}s non
lin{\'e}aire, pour diff{\'e}rentes valeurs du retard
r{\'e}guli{\`e}\-rement espac{\'e}es, et {\`a} choisir celle pour
lesquelles on approche le mieux la sortie. Cette m{\'e}thode est
inapplicable pour l'exemple $2$ avec la m{\'e}thode de moindre
carr{\'e} non lin{\'e}aire disponible dans scilab, car celle-ci ne
converge vers les param{\`e}tres que si elle est initialis{\'e}e avec
des valeurs voisines. Toutefois, les moindres carr{\'e}s
non-lin{\'e}aires permettent d'am{\'e}liorer la pr{\'e}cision de
l'{\'e}valuation des param{\`e}tres et du retard, une fois
initialis{\'e}s avec les r{\'e}sultats fournis par notre
m{\'e}thode. Le tableau ci-dessous fournit les valeurs approch{\'e}es
des {\'e}carts types obtenus par notre m{\'e}thode et apr{\`e}s
utilisation des moindres carr{\'e}s ($^{\ast}$), {\'e}valu{\'e}es sur
$100$ essais, la dur{\'e}e moyenne d'un calcul {\'e}tant de $10,5$~s.
\smallskip

\vbox{
{\scriptsize
\hbox to \hsize
{\hss\begin{tabular}{|c|c|c|c|c|c|c|c|c|c|c|c|c|}
\hline
\'E.~bruit & \'E.~$a_{1}$ & \'E.~$a_{1}$$^{\ast}$& \'E.~$a_{0}$ &
\'E.~$a_{0}$$^{\ast}$ & \'E.~$b$ & \'E.~$b$$^{\ast}$ & \'E.~$x(0)$ & \'E.~$x(0)$$^{\ast}$ & \'E.~$x'(0)$ & \'E.~$x'(0)$$^{\ast}$ &
\'E.~$h$ & \'E.~$h$$^{\ast}$\\
\hline
5 &  0,1   & 0,002 & 0,03  & 0,0005  & 0,2 & 0,002  & 1 & 0,4 & 6 & 0,5 & 0,07 & 0,001\\
\hline
10 & 0,15  & 0,005 & 0,06  & 0,001  & 0,35 & 0,005  & 1,7 & 0,6 & 10 & 0,9 & 0,1 & 0,0015\\
\hline
\end{tabular}\hss}
\hbox to \hsize{\hss Tableau 1. \'Ecart type des valeurs calcul{\'e}es\hss}
\hbox to \hsize{\hss \em Table 1. Standard deviations of computed values, stared values improved with mean-square\hss}}
}
\smallskip

Ce type de calcul suppose, naturellement l'identifiabilit{\'e} du
syst{\`e}me. Si la commande est nulle, la sortie constante etc. celui-ci
ne sera pas identifiable et il ne le sera que localement pour un
comportement p{\'e}riodique. Notons que, nos calculs reposant uniquement
sur des r{\'e}solutions au sens des moindres carr{\'e}s, il n'y a jamais
d'erreurs num{\'e}riques du type \og division par $0$.{\fg} Pour des retard
grands, les it{\'e}rations peuvent diverger et produire une erreur du type
{\tt invalid index} si les
intervales d'int{\'e}gration sortent de la p{\'e}riode d'observation du syst{\`e}me. 

\begin{exemple}
Cette m{\'e}thode se g{\'e}n{\'e}ralise au cas de plusieurs
retards. Consid{\'e}rons le syst{\`e}me $x''(t)+ax(t-h_{1})=bu(t-h_{2})$. En
utilisant $f(t)=1-\cos(t)$ sur $[0, 2\pi]$, on a $aI_{u,f,T_{1},T_{2}}$ 
$+ b\sin(h_2/(T_{2}-T_{1}))I_{x,\sin,T_{1},T_{2}}$ $+\cO(h_{2}^{2})$
$- bI_{x,f,T_{1},T_{2}}$
$- a\sin(h_1/(T_{2}-T_{1}))I_{x,\sin,T_{1},T_{2}}$ $+\cO(h_{1}^{2})$
$= (T_{2}-T_{1})^{-2}I_{x,f'',T_{1},T_{2}}$. En choisissant un nombre
suffisant d'intervalles $[T_{1},T_{2}]$, on peut {\'e}valuer par moindres
carr{\'e}s lin{\'e}aires des approximations de $a$,
$a\sin(h_1/(T_{2}-T_{1}))$, $b$ et $b\sin(h_1/(T_{2}-T_{1}))$. On en
d{\'e}duit des valeurs de $h_{1}$ et $h_{2}$ que l'on peut r{\'e}injecter dans
les {\'e}quations pour am{\'e}liorer la pr{\'e}cision.

Les simulations ont {\'e}t{\'e} faites avec $x=3\sin(t/2)+2\cos(t/3)$, $a=2,7$
$b=1,5$, $h_{1}=2$, $h_{2}=4$. L'echantillonage est effectu{\'e} {\`a}
$500~Hz$ avec des bruits gaussiens d'{\'e}cart-types diff{\'e}rents. On a pris
$T_{2}-T_{1}=10$ et $20$ intervalles avec $T_{1}$ variant par pas de
$2$ {\`a} partir de $T_{1}=15$. Le tableau suivant donne la moyenne et
l'{\'e}cart-type des valeurs calcul{\'e}es (sur 100 essais, chaque calcul
durant environ 1,5~s. avec 10 it{\'e}rations) pour diff{\'e}rents
niveaux de bruit. Les valeurs de $h_{1}$ et $h_{2}$ ne sont
naturellement d{\'e}finies que modulo la p{\'e}riode $12\pi$ de la sortie. La
pr{\'e}cision du calcul se d{\'e}grade si les intervalles n'ont pas une
longueur $T_{2}-T_{1}$ suffisament grande devant $h_{1}$ et $h_{2}$, mais aussi pour
un signal p{\'e}riodique si $T_{2}-T_{1}$ n'est pas assez petit devant la
p{\'e}riode.
\smallskip

\vbox{\em\scriptsize
\hbox to \hsize
{\hss\begin{tabular}{|c|c|c|c|c|c|c|c|c|}
\hline
\'E.~du bruit & M.~$a$ & \'E.~$a$ & M.~$b$ & \'E.~$b$ & M.~$h_{1}$ & \'E.~$h_{1}$ & M.~$h_{2}$ & \'E.~$h_{2}$\\
\hline
0,025 &  2,69 & 0,012 & 1,50 & 0,006 & 1,99 & 0,005 & 3,99 & 0,0051 \\
\hline
0,05  &  2,69 & 0,018 & 1,49 & 0,010 & 1,99 & 0,010 & 3,99 & 0,011\\
\hline
0,1   &  2,70 & 0,043 & 1,50 & 0,025 & 1,99 & 0,01  & 3,99 & 0,017\\
\hline
0,2   &  2,68 & 0,095 & 1,48 & 0,054 & 1,99 & 0,039 & 4,00 & 0,043\\
\hline
\end{tabular}\hss}
\hbox to \hsize{\hss Tableau 2. Moyenne et {\'e}cart type des valeurs calcul{\'e}es\hss}
\hbox to \hsize{\hss \em Table 2. Mean and standard deviations {\em (\'E.)} of computed values\hss}}

\end{exemple}

\section{Identification en temps r{\'e}el avec retard lentement variable}

On peut {\'e}galement adapter aux syst{\`e}mes {\`a} retard les m{\'e}thodes
d'identification et d'observation en temps r{\'e}el introduites dans la
section 1. On utilise les {\'e}quations
$$ \sum_{i=0}^{n} a_{i} \left(\sum_{k=j-1}^{n-1}
M_{i,j,k}I_{y,g_{k}}+\sum_{k=1}^{j} M_{i,n+j-1,k}I_{y,f_{k}}\right) +
b_{0}I_{u_{\hh},f_{j}} +b_{1}I_{u_{\hh},f_{j}} +0(h^{2})= 0,
$$ avec $b_{\ell}=bh^{\ell}$ et $u_{\hh}(t)=u(t-\hh)$. Ceci suppose que
$|h-\hh|\ll 1$. On peut s'y ramener quel que soit le retard par un
changement d'{\'e}chelle de temps. L'{\'e}valuation du retard est
alors $\hh+\hb_{1}/\hb_{0}$. On r{\'e}soud le syst{\`e}me par
int{\'e}gration au moyen d'une m{\'e}thode de gradient comme
d{\'e}crit section 2. Lorsque $t>T_{0}$, on suppose que les valeurs
des estimations sont devenues assez stables (il faut au moins que
$\hb_{0}$ ne s'annule plus). On pose alors $\hh'=0$ si $t<T_{0}$ et
$\hh'=\lambda_{h}b_{1}/b_{0}$ sinon. Pour plus de pr{\'e}cision, on
calcule $I_{u_{\hh},f_{j}}$ en posant
$I_{u_{\hh},f_{j}}'=(1+\lambda_{h}b_{1}/b_{0})(u_{\hh}(t)-\lambda
I_{u_{\hh},f_{j}})$.
\smallskip

\begin{exemple}
On utilise les m{\^e}mes commande, coefficients, type de bruit et
fr{\'e}quence que pour l'exemple 1. On prend $T_{0}=14$, $\lambda=1$,
$\Lambda=4.10^{-3}$ et $\lambda_{h}=0,4$.
\vskip 0.6truecm

\begin{figure}[!ht]\label{courbes-exemple3}
\begin{center}
\hbox to \hsize{\includegraphics[scale=0.35]{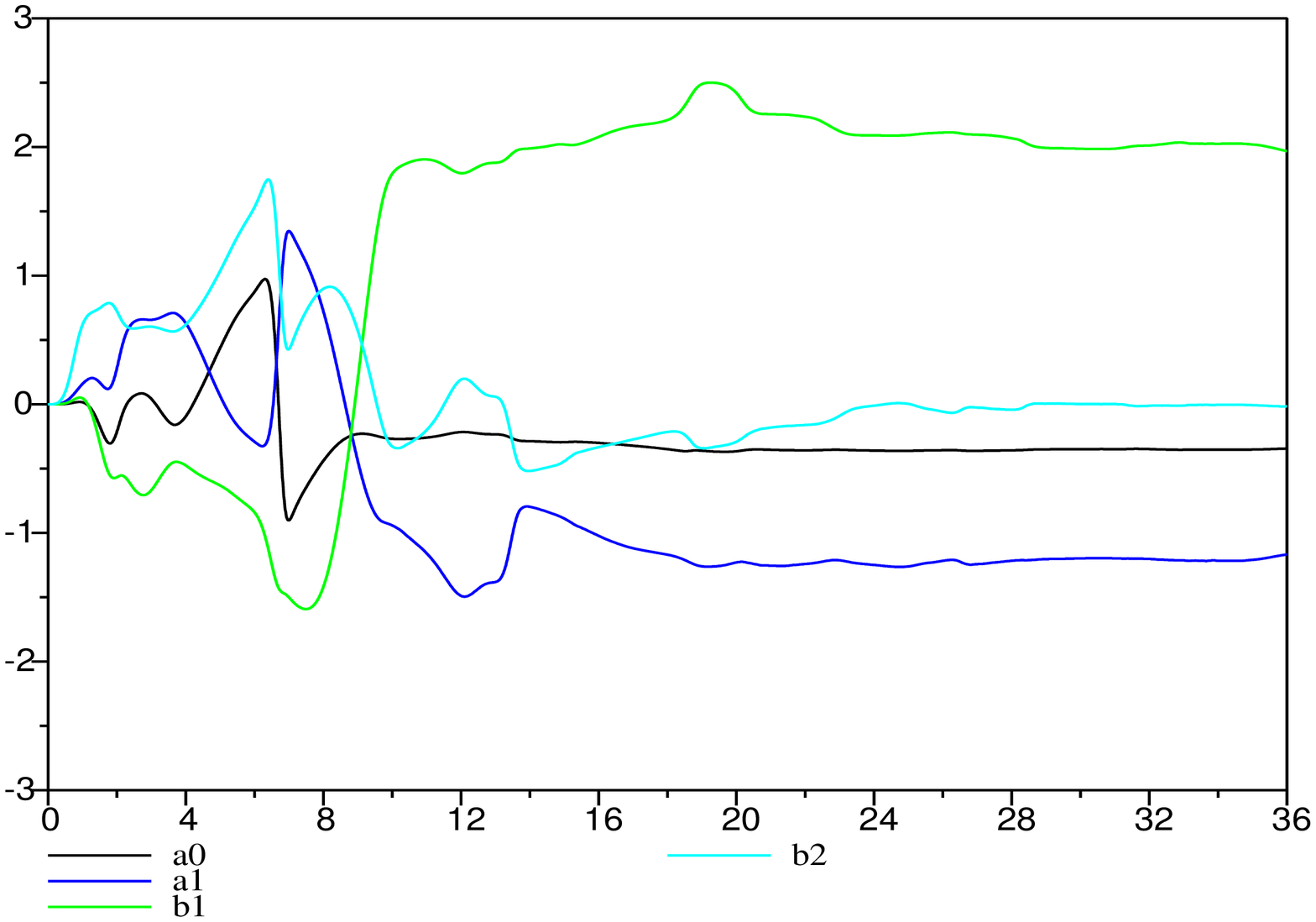}\hss
\includegraphics[scale=0.35]{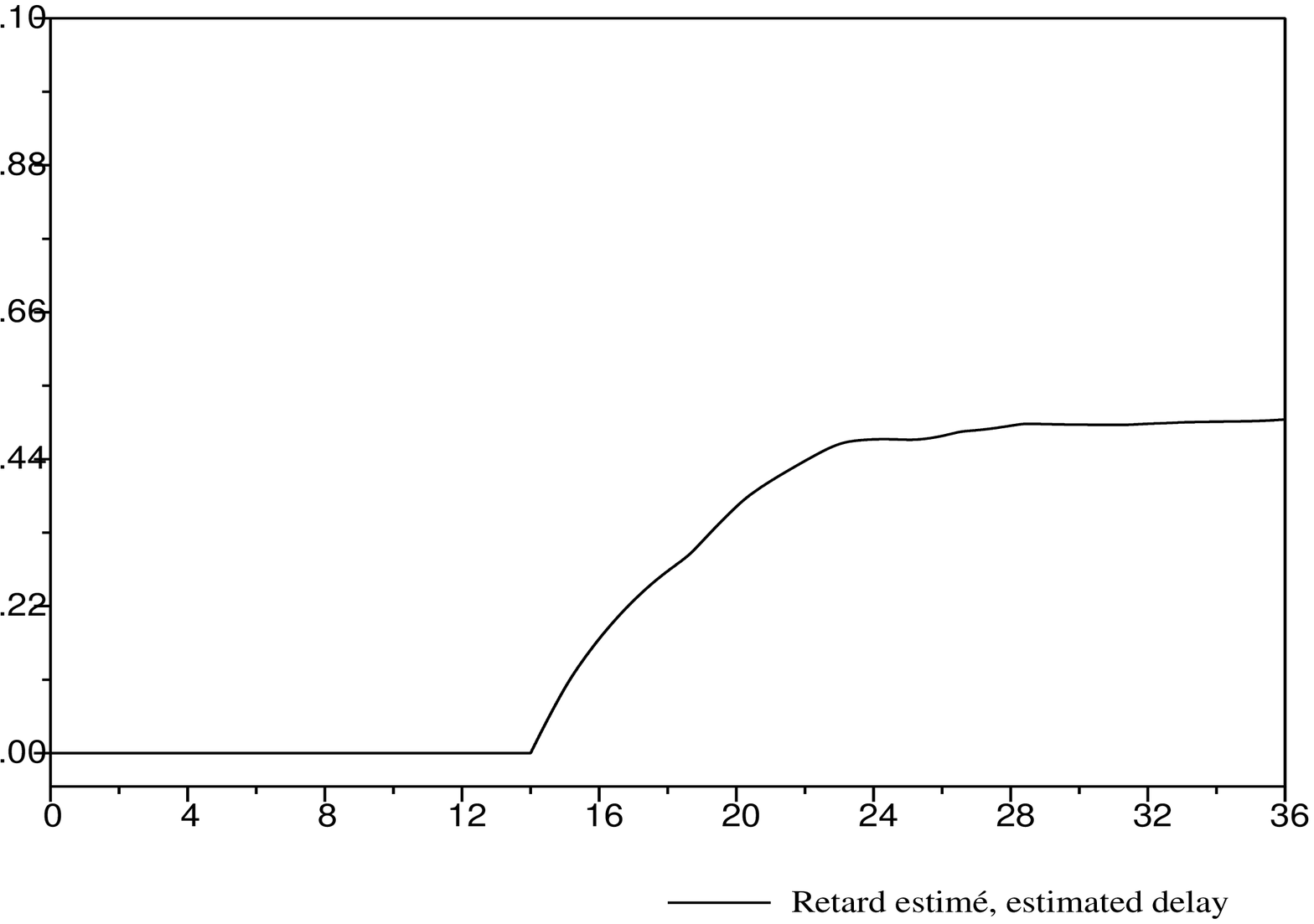}}
\caption{Estimations des coefficients et du retard, {\it Estimations
of coefficients and delay.}}
\end{center}
\end{figure}
\end{exemple}
\smallskip

\begin{exemple}
On utilise les m{\^e}mes commande, coefficients, type de bruit et
fr{\'e}quence que pour les exemples 1 et 3. On prend $T_{0}=40$,
$\lambda=1$, $\Lambda=4.10^{-4}$ et $\lambda_{h}=0,4$. Le retard est
initialis{\'e} {\`a} $3,6$ sa valeur {\`a} $T_{0}$, afin
d'acc{\'e}l{\'e}rer la convergence.
\begin{figure}[!ht]\label{courbes-exemple4}
\begin{center}
\hbox to \hsize{\includegraphics[scale=0.35]{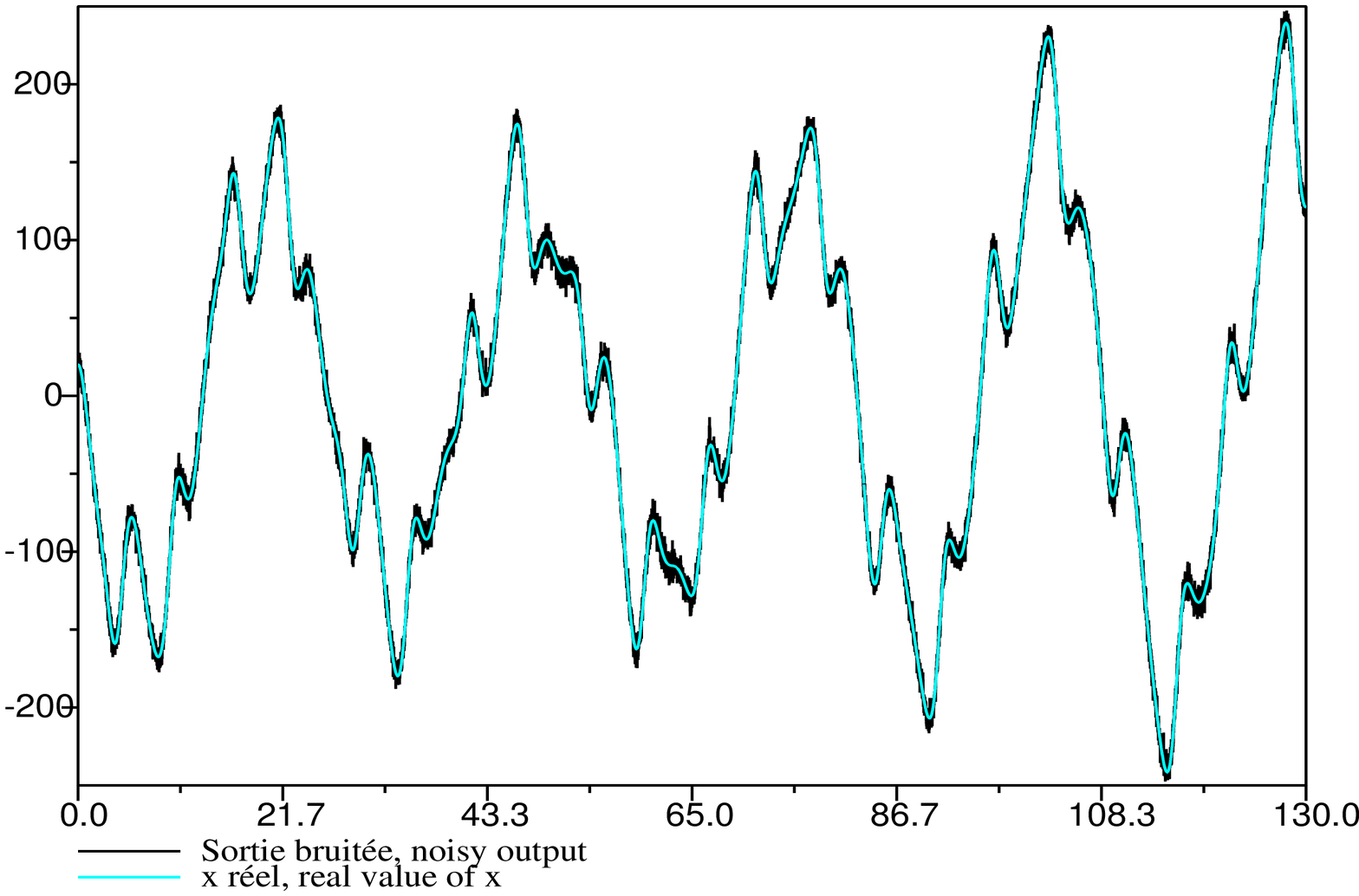}\hss
\includegraphics[scale=0.35]{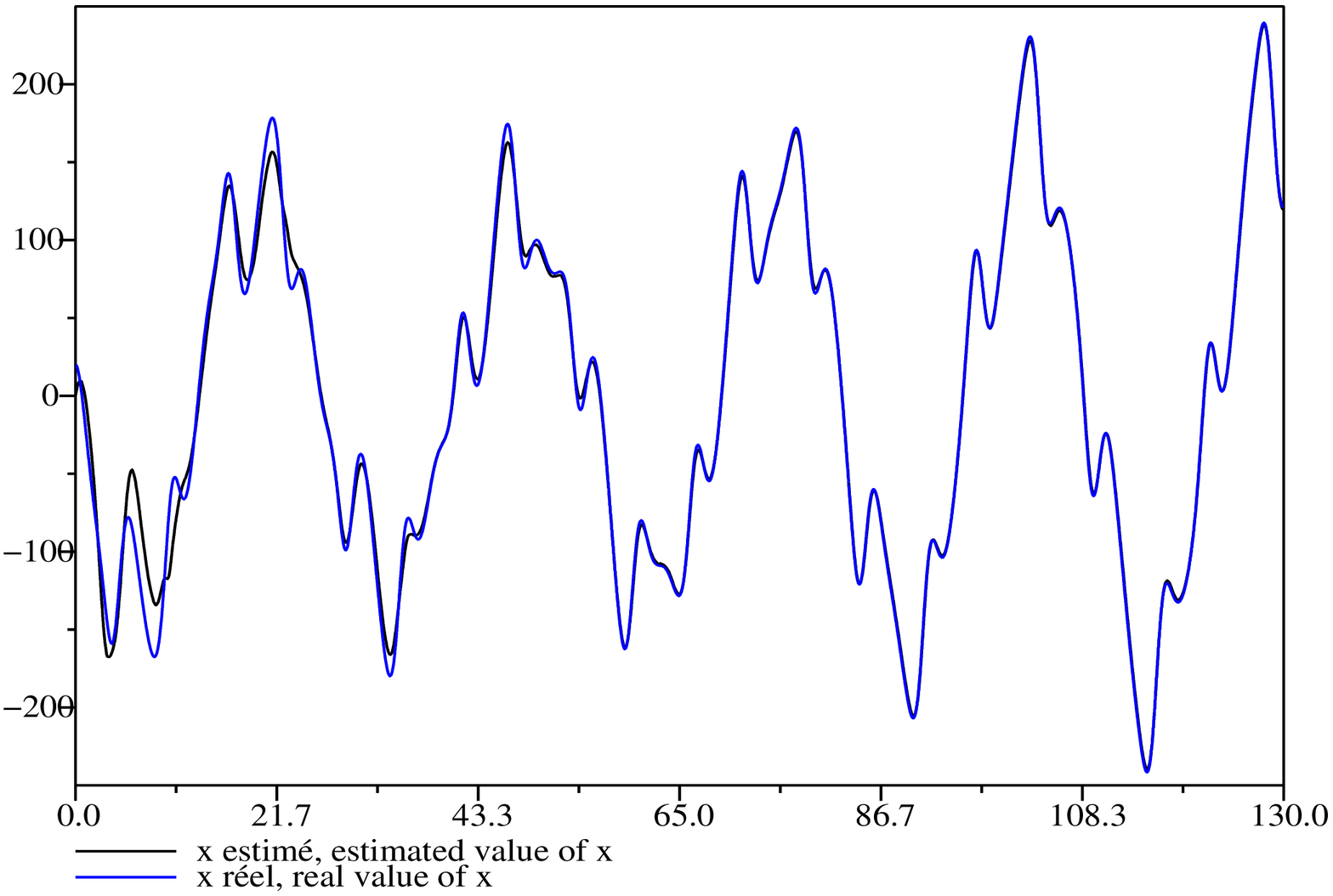}}
\medskip
\hbox to \hsize{\includegraphics[scale=0.35]{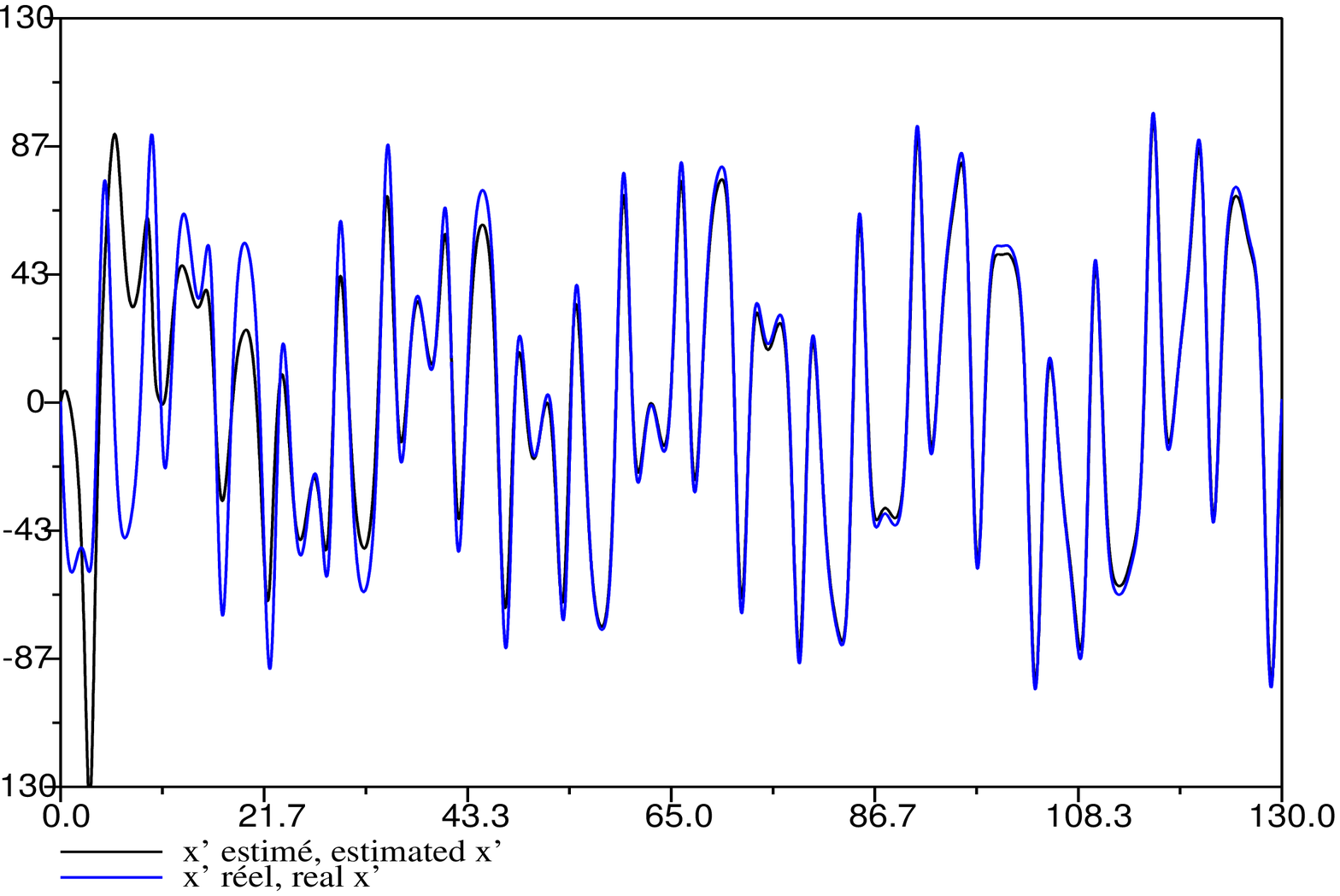}
\includegraphics[scale=0.35]{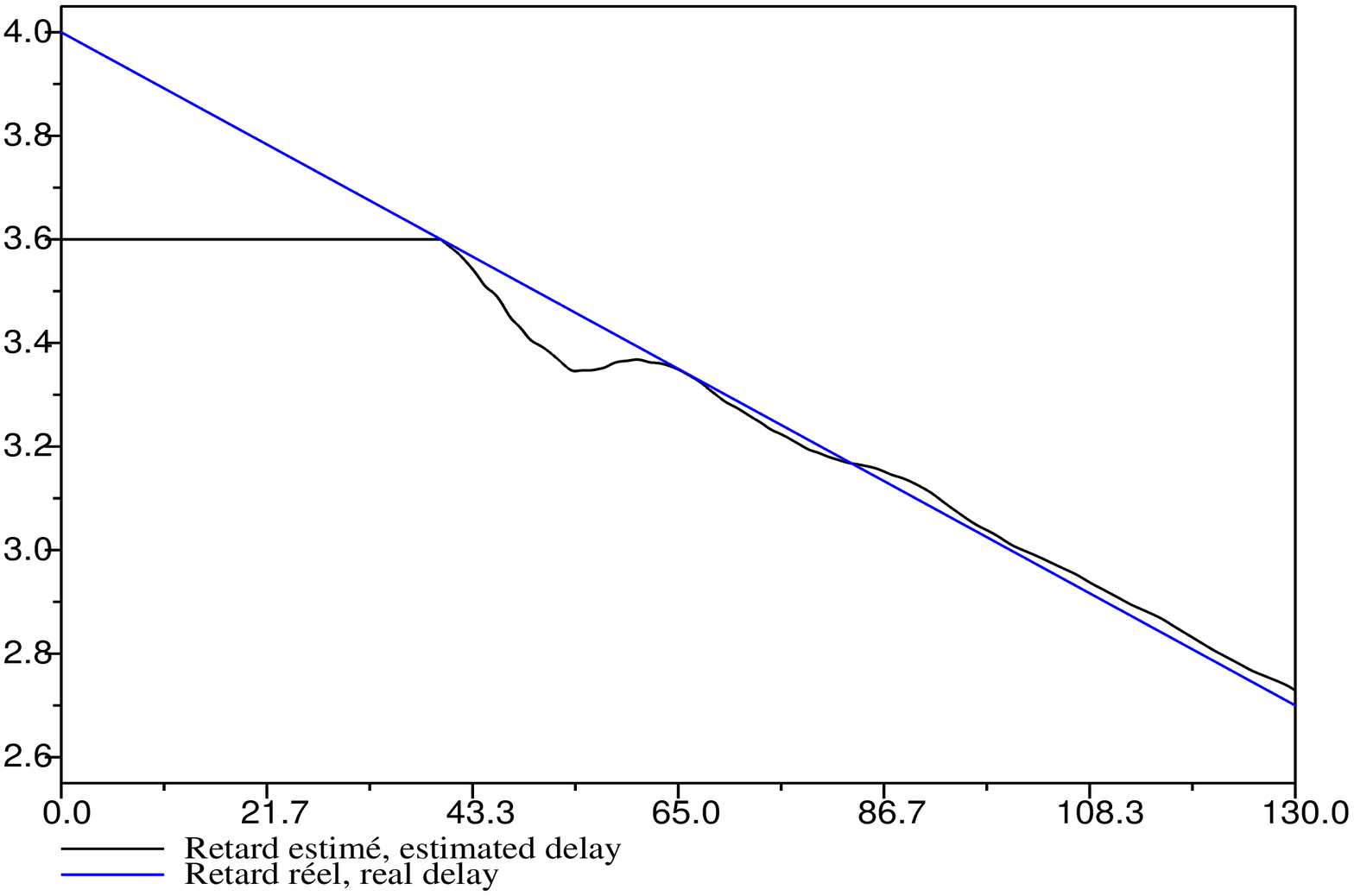}}
\caption{Estimations du retard de $x$ et de $x'$, {\it Estimations of the
delay, of $x$ and of $x'$.}}
\end{center}
\end{figure}
\end{exemple}

\section{Conclusion}

Les m{\'e}thodes que nous avons introduites permettent un large choix
dans leur mise en {\oe}uvre. Les quelques essais r{\'e}alis{\'e}s
montrent qu'elles sont applicables dans des situations r{\'e}alistes,
en particulier avec un niveau de bruit cons{\'e}quent. Elles peuvent
se g{\'e}n{\'e}raliser ais{\'e}ment, au moins en th{\'e}orie {\`a} des
situations plus complexes: retards multiples ou n'intervenant pas
uniquemement au niveau de la commande, comme l'illustre l'un de nos exemples.

Leur utilisation suppose la connaissance {\it a priori} d'une borne
sup{\'e}rieure pour le retard, ce qui est essentiel pour un choix appropri{\'e}
des diff{\'e}rents param{\`e}tres $\lambda$, $\Lambda$, $\lambda_{h}$
et $T_{0}$. En revanche, leur capacit{\'e} {\`a} s'adapter {\`a} des
coefficients ou {\`a} un retard lentement variable nous semble un atout
int{\'e}ressant. Les fichiers utilis{\'e}s pour les simulations sont
disponibles {\it via} l'url \cite{Simus}.



\begin{thebibliography}{00}

\bibitem{Belkoura-Richard-Fliess2006} Belkoura (L.), 
Richard (J.-P.), Fliess (M.) \og On-line identification of systems with
delayed inputs{\fg}, {\it MTNS'06, 16th Conf. Mathematical Theory of
Networks \& Systems}, 2006.



\bibitem{Diop-Kolmanovsky-Moraal-vanNieuwstadt2001} \textsc{Diop}
  (S.), \textsc{Kolmanovsky} (I.), \textsc{Moraal} (P.) et \textsc{van
  Nieuwstadt} (M.)  \og Preserving stability/performance when facing
  an unknown time delay{\fg}, {\it Control Engineering Practice}, {\bf
  9} 1319-1325, 2001.

\bibitem{Fliess-Mounier2000} \textsc{Fliess} (M.),
\textsc{Mounier} (H.), \og On a class of delay systems often arising in
practice{\fg}, {\it Kybernetika}, {\bf 27}, 295-308, 2000.

\bibitem{Fliess-Sira2003a} \textsc{Fliess} (M.) et
  \textsc{Sira-Ram{\`i}rez} (H.), \og An algebbraic  framework for
  linear identification{\fg}, {\it ESAIM Contr.\ Optim.\ Cal.\
  Variat.}, vol.~9, 151--168, 2003.

\bibitem{Fliess-Sira2003b} \textsc{Fliess} (M.) et
  \textsc{Sira-Ram{\`i}rez} (H.), \og Reconstructeurs d'{\'e}tat{\fg},
  {\it C.R.\ Acad.\ Sci.\ Paris}, 338 (2004), no. 1, 91--96.


\bibitem{Kolmanovskii-Myshkis1999} \textsc{Kolmanovskii} (V. B.)
et \textsc{Myshkis} (A.), {\it Introduction to the theory and
application of functionnal differential equations}, Dorctecht,
Kluwer Academy, 1999.

\bibitem{Mahji-Atherton1999} \textsc{Mahji} (S.) \textsc{Atherton}
  (D. P.), \og A novel identification method for time delay
  processes{\fg}, {\it ECC'99 (Fifth European Control Conference)},
  Karlsruhe, Allemangne, 1999.


\bibitem{Orlov-Belkoura-Richard-Dambrine2002} \textsc{Orlov} (Y.)
  \textsc{Belkoura} (L.) \textsc{Richard} (J. P.) et \textsc{Dambrine}
  (M.)  \og On-line parameter identification of linear time-delay
  systems{\fg}, {\it Proceedings of the 41$^{st}$ IEEE Conference on
  Decision and Control}, Las vegas, Nevada, USA, d{\'e}cembre 2002.

\bibitem{Richard2003}\textsc{Richard} (J.P.) \og Time-delay systems:
an overview of some recent advances and open problems{\fg} {\it
Automatica}, {\bf 39}, 10, 1667--1694, Octobre 2003.

\bibitem{Tan-Wang-Lee} \textsc{Tan} (K. K.) \textsc{Wang} (Q. K.) et
  \textsc{Lee} (T. II.), \og Finite spectrum assignment control of
  unstable time delay processes with relay tuning{\fg}, {\it
  Industrial Engineering and Chemical Research}, {\bf 37} (4),
  1351--1357, 1998.

\bibitem{Tuch-Feuer-Palmor1994} \textsc{Tuch} (J.) \textsc{Feuer} (A.)
  et \textsc{Palmor} (Z.J.), \og Time delay estimation in continuous
  linear time-invariant systems{\fg}, {\it IEEE Transactions on
  Automatic Control} {\bf 39}, 823--827.

\bibitem{Verduyn-Lunel1997} \textsc{Verduyn-Lunel} (S. M.), \og
  Identification problems in functional differential equations{\fg},
  {\it 36$^{\it th}$ IEEE CDC'99 (Conference on Decision and Control)},
  San Diego CA, d{\'e}cembre 1997, 4409--4413.


\bibitem{Simus} {\selectlanguage{english}\tt http://www.lix.polytechnique.fr/\~{ }ollivier/}

\end{thebibliography}
\end{document}